\documentclass[12pt]{article}

\usepackage{amssymb,amsmath,epsfig,graphics}

\newcounter{sscounter}[section]

\newtheorem{corollary}{Corollary}

\newtheorem{proposition}{Proposition}
\newtheorem{definition}{Definition}

\newtheorem{theorem}{Theorem}
\newtheorem{condition}{Condition}
\newcommand{\proof}{\noindent\textbf{Proof. }}
\newcommand{\qed}{{\hfill$\blacksquare$}}

\newcommand{\CC}{{\mathbb C }}
\newcommand{\GG}{{\mathbb G }}
\newcommand{\KK}{{\mathbb K }}
\newcommand{\LL}{{\mathbb L }}

\newcommand{\PP}{{\mathbb P }}

\newcommand{\RR}{{\mathbb R }}

\newcommand{\ZZ}{{\mathbb Z }}

\newcommand{\cA}{\mathcal{A}}

\newcommand{\cI}{\mathcal{I}}
\newcommand{\cL}{\mathcal{L}}

\newcommand{\cR}{\mathcal{R}}

\newcommand{\cX}{\mathcal{X}}
\newcommand{\cZ}{\mathcal{Z}}

\newcommand{\GL}{\mathcal{G}_\LL}
\newcommand{\oGL}{\overline{\mathcal{G}}_\LL}

\newcommand{\tGL}{\widetilde{\mathcal{G}}_{\LL}}
\newcommand{\CN}{\mathbb{C}^{N}}
\newcommand{\CNs}{{\CC^*}^{N}}

\newcommand{\PN}{\mathbb{P}^{N-1}}
\newcommand{\ZN}{\mathbb{Z}^{N}}

\newcommand{\KP}{\KK_\sP}

\newcommand{\sP}{\mathsf{P}}
\newcommand{\Pb}{\mathsf{P}^{2\bullet}}
\newcommand{\Pw}{\mathsf{P}^{2\circ}}
\newcommand{\Pn}{\mathsf{P}^{0}}
\newcommand{\Pe}{\mathsf{P}^{1}}

\newcommand{\va}{\mathbf{a}}
\newcommand{\vb}{\mathbf{b}}
\newcommand{\vc}{\mathbf{c}}

\newcommand{\ve}{\mathbf{e}}

\newcommand{\vu}{\mathbf{u}}

\newcommand{\vw}{\mathbf{w}}
\newcommand{\vx}{\mathbf{x}}
\newcommand{\vy}{\mathbf{y}}
\newcommand{\vz}{\mathbf{z}}
\newcommand{\vY}{\mathbf{Y}}
\newcommand{\vB}{\mathbf{B}}

\newcommand{\vl}{\mathfrak{l}}

\newcommand{\wt}{\varpi}

\newcommand{\gb}{\beta}

\newcommand{\ep}{\varepsilon}

\newcommand{\Hom}{\mathrm{Hom}}
\newcommand{\tors}{\mathrm{tors}}

\newcommand{\ol}[1]{\overline{#1}}

\newcommand{\modquot}[2]{\mbox{\raisebox{.2ex}{$#1$}\hspace{-.3em}/ \hspace{-.6em} \raisebox{-.2ex}{$#2$}}}

\begin{document}

\title{Chow Forms, Chow Quotients and\\ Quivers with Superpotential}

\author{Jan Stienstra\\
\small Mathematisch Instituut, Universiteit Utrecht, the Netherlands\\ 
\small e-mail: {J.Stienstra}{`at'}{uu.nl} \normalsize}

\date{}

\maketitle

\begin{abstract}
We recast the correspondence between 3-dimensional toric Calabi-Yau singularities and quivers with superpotential in the setting of an abelian algebraic group $\GL$ acting on a linear space $\CN$.
We show how the quiver with superpotential gives a simple explicit description of the Chow forms of the closures of the orbits
in the projective space $\PN$.
The resulting model of the orbit space is known as 
the Chow quotient of $\PN$ by $\GL$. 
As a by-product we prove a result about the relation between
the quiver with superpotential
and the principal $\cA$-determinant in Gelfand-Kapranov-Zelevinsky's theory
of discriminants.
\end{abstract}

\section{Introduction.}\label{section:introduction}
One of the important themes in the physics literature on AdS/CFT is a correspondence between complex $3$-dimensional Calabi-Yau singularities
and quiver gauge theories. A breakthrough in this field was the discovery
of an algorithm that realizes the correspondence for toric CY3 singularities
\cite{HV,Ky}. Experience with, for instance, the Gelfand-Kapranov-Zelevinsky theory of hypergeometric systems or De Bruijn's construction of Penrose tilings
shows that sometimes great simplifications and new insights can be achieved by passing to an equivariant setting in higher dimensions. 
In \cite{S} we adapted the aforementioned algorithm to a higher dimensional viewpoint.
In the present paper we clarify some equivariance aspects.
We study an $(N-2)$-dimensional algebraic subgroup 
$\GL$ of $\CNs$ acting 
by diagonal matrices on $\CN$ and the induced action of a
subgroup $\tGL$ of $\GL$ on $\CN$ and of a
quotientgroup $\oGL$ of $\GL$ 
on the projective space $\PN$. Both groups $\tGL$ and $\oGL$ have dimension
$N-3$.
We make models of the orbit spaces by embedding
$\modquot{\CNs}{\tGL}$ into some linear space and $\modquot{\CNs}{\GL}$
into some projective space and taking closures.
We construct these embeddings directly
from the quiver  with superpotential $\sP$ on the other side of the correspondence.
For this we reformulate the data of $\sP$
in what we call the bi-adjacency matrix $\KP(\vz,\vu)$; see Definition \ref{def:bi-adjacency matrix}. Here $\vz$ denotes a tuple of variables with for every quiver arrow $e$ one variable $z_e$
and $\vu$ denotes a tuple of variables with for every quiver node $i$ one variable $u_i$. The entries of $\KP(\vz,\vu)$ are elements in the polynomial ring $\ZZ[\vz,\vu]$.
Viewing $u_1,\ldots,u_N$ as coordinates on $\CN$ we get an action of 
the group $\GL$ on the polynomial ring $\CC[\vu]=\CC[u_1,\ldots,u_N]$. 
We show in Theorem \ref{thm:homogeneity} that under suitable conditions on $\sP$ and $\GL$ there is a character $\chi$ of $\GL$ such that
in the polynomial ring $\CC[\vz,\vu]$
$$
\det\KP(\vz,\xi\vu)=\chi(\xi)\det\KP(\vz,\vu)\,,\qquad\forall \xi\in\GL\,. 
$$
 The above mentioned group $\tGL$ is the kernel of $\chi$. By evaluating $\vu=(u_1,\ldots,u_N)$ at the points of
$\CNs$ we obtain well defined maps
\begin{eqnarray}
\label{eq:quot1}
\modquot{\CNs}{\tGL}&\longrightarrow&\CC[\vz]^{(\nu)}\,,
\qquad\quad\: \vu\mapsto\det\KP(\vz,\vu)\,,
\\
\label{eq:quot2}
\modquot{\CNs}{\GL}&\longrightarrow&\PP(\CC[\vz]^{(\nu)})\,,
\qquad \vu\mapsto\det\KP(\vz,\vu)\bmod\CC^*\,;
\end{eqnarray}
here $\CC[\vz]^{(\nu)}$ is the homogeneous part of
degree $\nu=\deg_\vz\det\KP(\vz,\vu)$ in the polynomial ring $\CC[\vz]$
and $\PP(\CC[\vz]^{(\nu)})$ is the corresponding projective space.
Taking Zariski closures of the images of the maps (\ref{eq:quot1})
and (\ref{eq:quot2}) one obtains models for the orbit spaces of
$\tGL$ acting on $\CN$ and $\oGL$ acting on $\PN$ together with embeddings into $\CC[\vz]^{(\nu)}$ and $\PP(\CC[\vz]^{(\nu)})$, respectively.

In Definition \ref{def:complementary} we define
the \emph{complementary bi-adjacency matrix} $\KP^c(\vz,\vu)$
and in Formula (\ref{eq:pluckerweight}) we define a ring homomorphism
$\vy$ from the polynomial ring $\ZZ[\vz]$ to  
the ring $\cR_{2,N}$ of homogeneous coordinates on the Grassmannian
$\GG(2,N)$.
In Theorem \ref{theorem:equal} we show that for fixed $\vu\in\CNs$
the element $\det\KP^c(\vy(\vz),\vu)$ of $\cR_{2,N}$
is a \emph{Chow form} for the orbit closure $\ol{\GL[\vu]}$ in $\PN$.
The maps (\ref{eq:quot1}) and (\ref{eq:quot2}) have obvious analogues:
\begin{eqnarray}
\label{eq:quot3}
\modquot{\CNs}{\tGL}&\longrightarrow&\cR_{2,N}^{(\nu)}\,,
\qquad\quad\: \vu\mapsto\det\KP^c(\vy(\vz),\vu)\,,
\\
\label{eq:quot4}
\modquot{\CNs}{\GL}&\longrightarrow&\PP(\cR_{2,N}^{(\nu)})\,,
\qquad \vu\mapsto\det\KP^c(\vy(\vz),\vu)\bmod\CC^*\,;
\end{eqnarray}
The closure of the image of the map (\ref{eq:quot4}) in 
$\PP(\cR_{2,N}^{(\nu)})$ is a model for the orbit space of
$\oGL$ acting on $\PN$ known as the \emph{Chow quotient}; see \cite{KSZ}.
In \cite{KSZ} one can also find a discussion of how this model compares with GIT quotients.

A priori $\det\KP^c(\vy(\vz),\vu)$ is an element in 
$\cR_{2,N}[\vu]=\cR_{2,N}[u_1,\ldots,u_N]$, i.e.
in the ring of homogeneous coordinates on $\GG(2,N)\times\PN$.
In (\ref{eq:quot4}) this element is evaluated at points of $\PN$.
One can also evaluate it at points of $\GG(2,N)$
and find polynomials in the ring $\CC[\vu]=\CC[u_1,\ldots,u_N]$.
A point of $\GG(2,N)$ is a line in $\PN$ and the corresponding polynomial is the equation for the image of this line in the Chow quotient for
$\oGL$ acting on $\PN$.   

The $\LL$ which appears in the above presentation as a subscript is in fact the character lattice of the group $\modquot{\CNs}{\GL}$. It is a rank $2$ subgroup of $\ZN$ and thus defines a point $\vl$ in the Grassmannian
$\GG(2,N)$. 
In Theorem \ref{thm:Adeterminant} we show that evaluating 
$\det\KP^c(\vy(\vz),\vu)$ at the point
$\vl$ yields the principal $\cA$-determinant of Gelfand-Kapranov-Zelevinsky \cite{gkz4}, if $\GL$ is connected.
In \cite{S} this result was observed in some examples and then conjectured to hold in general.

In Figures \ref{fig:fan F1}, \ref{fig:dP3I}, \ref{fig:bi-ad}
and Section \ref{section:example} we give one example of the theory
presented in this paper. Some further examples can be found in \cite{S}.

\section{The lattice $\LL$ and associated groups.}\label{section:groups}
Throughout this paper $\LL$ is a rank $2$ subgroup of $\ZN$ not contained in any of the standard coordinate hyperplanes and contained in the kernel of the map
\begin{equation}\label{eq:h}
h:\ZN\longrightarrow\ZZ\,,\qquad h(z_1,\ldots,z_N)=z_1+\ldots+z_N\,.
\end{equation}
Unlike what seems to be practice in the theory of Gelfand, Kapranov and Zelevinsky \cite{gkz4}, we do allow torsion in the quotient group
$\modquot{\ZN}{\LL}$ .
Denoting by $\ve_1,\ldots,\ve_N$ the standard basis of
$\ZZ^N$ we set
$$
\GL:=\Hom(\modquot{\ZN}{\LL},\CC^*)\,,\qquad
\va_i:= \; \ve_i \bmod\LL\;
\textrm{ in }\;\modquot{\ZN}{\LL}\,.
$$
$\GL$ is an algebraic group for which the connected component 
containing the
identity is isomorphic to the complex torus $(\CC^*)^{N-2}$ and the group of connected components is isomorphic to the
finite abelian group $(\modquot{\ZN}{\LL})_\tors$.
The inclusion $\LL\hookrightarrow\ZN$ induces the inclusion of groups
$$
\GL\hookrightarrow \CNs\,,\qquad\xi\mapsto (\xi(\va_1),\xi(\va_2),\ldots,\xi(\va_N))
$$
and thus an action of $\GL$ on $\CN$ and $\PN$.  
The map $h$ induces a homomorphism 
$\ol{h}:\modquot{\ZN}{\LL}\rightarrow\ZZ$ such that
$\ol{h}(\va_i)=1$ for all $i$. As a consequence there is a subgroup
$\GL'$ of $\GL$ consisting of those
$\xi\in\GL$ for which $\xi(\va_1)=\ldots=\xi(\va_N)$.
So, in fact $\GL$ acts on $\PN$ via the quotient group 
$\oGL:=\modquot{\GL\,}{\GL'}$.

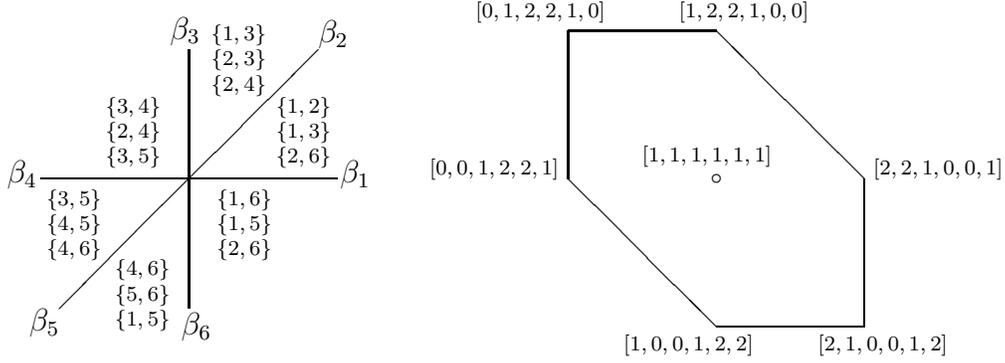
\begin{figure}[t]
\setlength{\unitlength}{0.7pt}
\begin{picture}(320,160)(-130,-95)
\put(-25,0){
\begin{picture}(320,160)(0,0)
\put(0,0){\line(1,0){80}}
\put(0,0){\line(-1,0){80}}
\put(0,0){\line(0,1){70}}
\put(0,0){\line(0,-1){70}}
\put(0,0){\line(1,1){70}}
\put(0,0){\line(-1,-1){70}}
\put(85,-3){\makebox(10,10){$\gb_1$}}
\put(-95,-3){\makebox(10,10){$\gb_4$}}
\put(-7,75){\makebox(10,10){$\gb_3$}}
\put(-1,-85){\makebox(10,10){$\gb_6$}}
\put(73,73){\makebox(10,10){$\gb_2$}}
\put(-83,-83){\makebox(10,10){$\gb_5$}}
\scriptsize
\put(-50,10){\makebox(40,30){$\begin{array}{l}
\{3,4\}\\ \{2,4\}\\ \{3,5\}\end{array}$}}
\put(10,-40){\makebox(40,30){$\begin{array}{l}
\{1,6\}\\ \{1,5\}\\ \{2,6\}\end{array}$}}
\put(42,10){\makebox(40,30){$\begin{array}{l}
\{1,2\}\\ \{1,3\}\\ \{2,6\}\end{array}$}}
\put(-82,-40){\makebox(40,30){$\begin{array}{l}
\{3,5\}\\ \{4,5\}\\ \{4,6\}\end{array}$}}
\put(-45,-78){\makebox(40,30){$\begin{array}{l}
\{4,6\}\\ \{5,6\}\\ \{1,5\}\end{array}$}}
\put(7,49){\makebox(40,30){$\begin{array}{l}
\{1,3\}\\ \{2,3\}\\ \{2,4\}\end{array}$}}
\end{picture}
}
\put(260,0){
\begin{picture}(320,160)(0,0)
\put(80,0){\line(-1,1){80}}
\put(-80,0){\line(1,-1){80}}
\put(80,0){\line(0,-1){80}}
\put(-80,0){\line(0,1){80}}
\put(0,80){\line(-1,0){80}}
\put(0,-80){\line(1,0){80}}
\put(0,0){\circle{5}}
\scriptsize
\put(115,0){\makebox(10,10){$[2,2,1,0,0,1]$}}
\put(-125,0){\makebox(10,10){$[0,0,1,2,2,1]$}}
\put(10,85){\makebox(10,10){$[1,2,2,1,0,0]$}}
\put(-20,-95){\makebox(10,10){$[1,0,0,1,2,2]$}}
\put(-100,85){\makebox(10,10){$[0,1,2,2,1,0]$}}
\put(85,-95){\makebox(10,10){$[2,1,0,0,1,2]$}}
\put(-10,7){\makebox(10,10){$[1,1,1,1,1,1]$}}
\end{picture}
}
\end{picture}
\caption{\label{fig:fan F1}\textit{
Secondary fan with lists $L_\vc$ (left)
and secondary polygon with vectors $\sum_{\{i,j\}\in L_\vc}(\ve_i+\ve_j)$
(right)
for 
\mbox{$\LL=\ZZ^2$ \footnotesize{$\left[\protect\begin{array}{rrrrrr}
\!1&\!1&\!0&\!-1&\!-1&\!0\protect\\ \!0&\!1&\!1&\!0&\!-1&\!-1\protect\end{array}\right]$}}.
\normalsize 
\protect\\
This shows $\va_0=(1,1,1,1,1,1)\bmod\LL$.}}
\setlength{\unitlength}{1pt}
\normalsize
\end{figure}

The subgroup $\tGL$ of $\GL$ is the kernel of some character $\chi$ of $\GL$. The construction is slightly involved.
We choose a basis for $\LL$ and represent the inclusion $\LL\hookrightarrow\ZN$ by a matrix 
$B=(b_{ij})_{i=1,2;j=1,\ldots,N}$ with entries in $\ZZ$.
Let $\gb_1,\ldots,\gb_N\in\ZZ^2$ denote the columns of $B$.
The complete fan in $\RR^2$ with $1$-dimensional rays
$\RR_{\geq 0}\gb_i$ ($i=1,\ldots,N$) is the so-called \emph{secondary fan} associated with $\LL$;
see Figure \ref{fig:fan F1} for an example.
Take $\vc\in\RR^2\setminus\bigcup_{i=1}^N\RR_{\geq 0}\gb_i$ and set
$L_\vc\,=\,\{\,\{i,j\}\subset\{1,\ldots,N\}\:|\:\vc\in\RR_{\geq 0}\gb_i+
\RR_{\geq 0}\gb_j\:\}$.
We now define:
\begin{eqnarray}
\label{eq:a0}
\va_0&:=&\textstyle{\sum_{\{i,j\}\in 
L_\vc}|\det(\gb_i,\gb_j)|(\va_i+\va_j)}\,,
\\
\label{eq:chi}
\chi&:&\GL\rightarrow\CC^*\,,\quad \chi(\xi):=\xi(\va_0)\,,
\\
\label{eq:tGL}
\tGL&:=&\ker\chi\,.
\end{eqnarray}
It follows from \cite{S} Eqs.(6), (18) that $\va_0$ does 
not depend on the particular choice of the vector $\vc$.
Changing the basis for $\LL$ multiplies $B$ with a $2\times 2$-matrix with determinant $\pm 1$. So, $\va_0$ 
does also not depend on the particular choice of the basis for $\LL$.

\section{The quiver with superpotential $\sP$.}\label{section:quiver}
On the quiver side we look at a quiver with a special kind of superpotential $\sP=(\Pn,\Pe,\Pb,\Pw, s, t, \vb, \vw)$.
It consists of four finite sets $\sP^0,\sP^1,\Pb,\Pw$ 
and four maps
\begin{equation}\label{eq:P}
s,\,t\,:\:\Pe\rightarrow\Pn\,,\qquad \vb:\Pe\rightarrow\Pb
\,,\qquad \vw:\Pe\rightarrow\Pw\,.
\end{equation}
The elements of $\Pn$ are called $0$-cells or nodes or vertices.
The elements of $\Pe$ are called $1$-cells or arrows or edges. 
The elements of $\Pb$ (resp. $\Pw$) are called black (resp. white) $2$-cells. 
An arrow $e$ is oriented from $s(e)$ to $t(e)$.

\emph{Througout this paper the following condition is assumed to be satisfied:}
\begin{condition}\label{condition P}.

\begin{enumerate}
\item 
The directed graph $(\Pn,\Pe,s,t)$ is connected and has no oriented cycles of length $\leq 2$ and in every node $v\in\Pn$ there are as many incoming as outgoing arrows; i.e.
$\sharp\{e\in\Pe|t(e)=v\}\,=\,\sharp\{e\in\Pe|s(e)=v\}$.
\item 
There are as many black as white cells; i.e. $\sharp\Pb=\sharp\Pw$.
\item 
For every $\vb\in\Pb$ and every $\vw\in\Pw$
the sets $\{e\in\Pe\,|\,\vb(e)=\vb\}$ and $\{e\in\Pe\,|\,\vw(e)=\vw\}$
are connected oriented cycles; by this we mean that the elements can be ordered
$(e_1,e_2,\ldots,e_r)$ such that
$t(e_i)=s(e_{i+1})$ for $i=1,\ldots,r-1$ and $t(e_r)=s(e_1)$.
\end{enumerate}
\end{condition}

One can realize $\sP$ geometrically as an oriented surface without boundary by
representing every black $2$-cell $\vb$ (resp. white $2$-cell $\vw$)
as a convex polygon with sides labeled by the elements $e\in\Pe$ 
with $\vb(e)=\vb$ (resp. $\vw(e)=\vw$) and by gluing these polygons along 
sides with the same label.
In quiver theory one views $(\Pn,\Pe,s,t)$ as a quiver (=directed graph) and the
boundary cycles of the $2$-cells as the terms of a superpotential; cf. \cite{S} 8.6.

\begin{figure}[t]
\begin{picture}(350,100)(-40,-10)
\put(0,0){
\begin{picture}(150,100)(0,15)
\setlength\epsfxsize{3.5cm}
\epsfbox{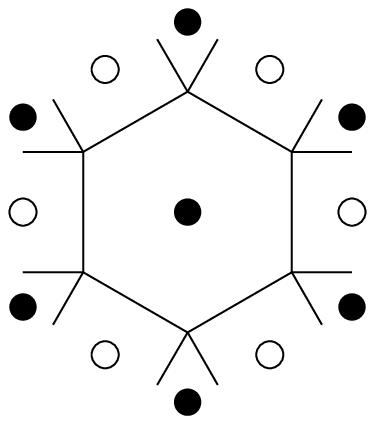}
\end{picture}
}

\put(180,-10){
\begin{picture}(150,100)(0,15)
\setlength\epsfxsize{4cm}
\epsfbox{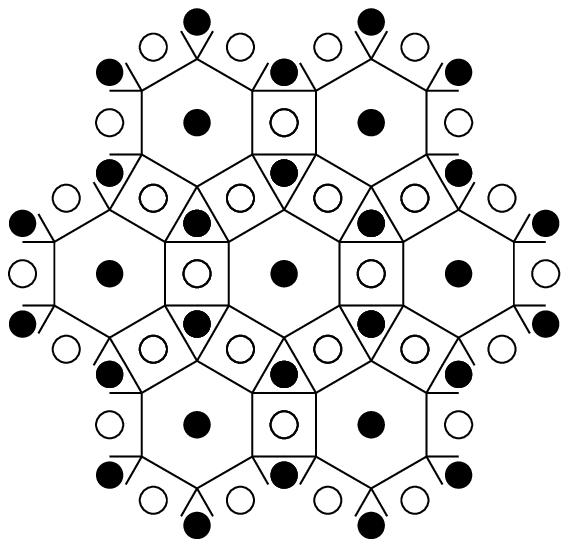}
\end{picture}
}
\end{picture}
\caption{\label{fig:dP3I}
\textit{On the left: the quiver with superpotential $\sP$; 
opposite sides of the hexagon must be identified; $\sP$ actually lies on a torus. On the right: part of the lift of $\sP$ on the universal cover of the torus.}}
\end{figure}

The structure of $\sP$ described in (\ref{eq:P}) can be accurately captured in
its bi-adjacency matrix:

\begin{definition}\label{def:bi-adjacency matrix}
The 
\emph{bi-adjacency matrix} of $\sP$
is the square matrix $\KP(\vz,\vu)$ with rows (resp. columns)
corresponding with the black (resp. white) $2$-cells,
with entries in the polynomial ring
$\ZZ[z_e,u_i\:|\:e\in\Pe\,,\;i\in\Pn\,]$ and
$(\vb,\vw)$-entry: 
\begin{equation}\label{eq:bi-adjacency matrix}
(\KP(\vz,\vu))_{\vb,\vw}\,=\,\sum_{e\in \sP^1:\,\vb(e)=\vb,\,\vw(e)=\vw}
z_e\,u_{s(e)}u_{t(e)}\,.
\end{equation}
\end{definition}

This definition of the bi-adjacency matrix is the same as \cite{S} Eq.(1), except that in loc. cit. we used a $\CC$-valued function $\wt$ on $\Pe$ instead of the formal variables $z_e$.
The products $u_{s(e)}u_{t(e)}$ of commuting variables do not show the orientation of the edge $e$. As explained in \cite{S} \S 8.10, the orientation can easily be reconstructed. So there is no loss of information in going from
(\ref{eq:P}) to (\ref{eq:bi-adjacency matrix}).
The bi-adjacency matrix $\KP(\vz,\vu)$ of $\sP$ is also the \emph{Kasteleyn matrix of a twist
of $\sP$}; see \cite{S} Section 8.

\section{Compatibility of $\LL$ and $\sP$.}\label{section:compatibility}
In \cite{S} \S6 we described an algorithm, based on 
the \emph{Fast Inverse Algorithm} of
\cite{HV} and the \emph{untwisting procedure} of \cite{FHKV},
for constructing from a rank $2$ lattice $\LL$ as in Section \ref{section:groups} 
a quiver with superpotential $\sP$ as in Section \ref{section:quiver},
such that also Condition \ref{condition LP} (below) is satisfied.
We isolated Conditions \ref{condition P} and \ref{condition LP} 
because this is exactly what we need in the present paper
and because these conditions can be checked without worrying about 
the algorithm of \cite{S}.

\begin{condition}\label{condition LP}
For every $\vb\in\Pb$ and $\vw\in\Pw$
there are elements $\ep_\vb$ and $\ep_\vw$ in $\modquot{\ZN}{\LL}$ such that
for every $e\in\sP^1$ and with $\va_0$ as in (\ref{eq:a0}) one has 
\begin{equation}\label{eq:cycle condition}
\va_{s(e)}\,+\,\va_{t(e)}=\ep_{\vb(e)}\,-\,\ep_{\vw(e)}\,,
\qquad
\va_0=\sum_{\vb\in\Pb}\ep_\vb\,-\,\sum_{\vw\in\Pw}\ep_\vw\,.
\end{equation}
\end{condition}

The ring $\CC[u_i\:|\:i\in\Pn]\,=\,\CC[u_1,\ldots,u_N]$
is the coordinate ring of $\CC^N$ and therefore carries an action of the group $\GL$:
$$
\xi\,u_i\,:=\,\xi(\va_i)u_i\qquad\textrm{for}\quad i=1,\ldots,N,\;\xi\in\GL\,.
$$

\begin{theorem}\label{thm:homogeneity}
Assume $\LL$ and $\sP$ satisfy Conditions \ref{condition P} and \ref{condition LP}. 
Then we have for every $\xi\in\GL$:
\begin{eqnarray}
\label{eq:bi-ad homogeneity}
\KP(\vz,\xi\vu)&=&\mathrm{diag(\xi(\ep_\vb))}\:\KP(\vz,\vu)\:
\mathrm{diag(\xi(-\ep_\vw))}\,,
\\
\label{eq:det homogeneity}
\det\KP(\vz,\xi\vu)&=&\chi(\xi)\det\KP(\vz,\vu)\,,
\\
\label{eq:sizebw}
\deg_\vz\det\KP(\vz,\vu)&=&
\textstyle{\frac{1}{2}} \deg_\vu\det\KP(\vz,\vu)=
\sharp\Pb=\sharp\Pw=\textstyle{\frac{1}{2}}\ol{h}(\va_0)\,.
\end{eqnarray}
with character
$\chi$ as in (\ref{eq:chi})and the entries of the diagonal matrices $\mathrm{diag(\xi(\ep_\vb))}$
and $\mathrm{diag(\xi(-\ep_\vw))}$
labeled with the black and white $2$-cells in agreement with the labeling of the rows and columns of $\KP(\vz,\vu)$. 
\end{theorem}
\proof (\ref{eq:bi-ad homogeneity}), (\ref{eq:det homogeneity}) 
and the first three equalities in (\ref{eq:sizebw}) are obvious.
The first half of (\ref{eq:cycle condition}) implies that there is a $k\in\ZZ$ such that $\ol{h}(\ep_\vb)=\ol{h}(\ep_\vw)+2=k$
for all $\vb\in\Pb$ and $\vw\in\Pw$. The last equality in
(\ref{eq:sizebw}) now follows from the second half of (\ref{eq:cycle condition}).
\qed

\begin{figure}[t]
\begin{picture}(350,70)(-10,-20)
\small
\put(0,0){
\mbox{
\begin{tabular}{r||c|c|c|}
$\ep_\vb\downarrow\,:\,\ep_\vw\rightarrow$&$(0,-1,-1,0,0,0)$&
$(-1,-1,0,0,0,0)$&$(0,0,-1,-1,0,0)$\\ \hline\hline
$(0,0,0,0,0,0)$&$z_1u_2u_3+z_2u_5u_6$&$z_3u_1u_2+z_4u_4u_5$&$
z_5u_3u_4+z_6u_1u_6$ \\ \hline
$(0,0,-1,0,0,1)$&$z_7u_2u_6$&$z_8u_2u_4$&$z_9u_4u_6$ \\ \hline 
$(0,-1,0,0,1,0)$&$z_{10}u_3u_5$&$z_{11}u_1u_5$&$z_{12}u_1u_3$ \\ \hline 
\end{tabular}
}}
\end{picture}
\caption{\label{fig:bi-ad}
\textit{The bi-adjacency matrix for $\sP$ as in Figure \ref{fig:dP3I}
and on its sides vectors (to be read modulo the $\LL$ of Figure 
\ref{fig:fan F1}) which verify Condition \ref{condition LP}.}}
\end{figure}

\section{Lines in projective space and Chow forms.}\label{section:lines}
Every line in $\PP^{N-1}$ is of the form
$$
\cL_Y:=\{[u_1:\ldots:u_N]\in\PP^{N-1}\:|\: u_j=
t_1y_{1j}+t_2y_{2j}\,,\; [t_1:t_2]\in\PP^1\;\}
$$
with $Y=(y_{ij})_{i=1,2;j=1,\ldots,N}$ a $2\times N$ complex matrix of rank $2$.
Let $M(2,N)$ denote the set of $2\times N$ complex matrices of rank $2$.
As $\cL_{Y'}=\cL_Y$ if and only if there is an
invertible $2\times 2$-matrix $g$ such that $Y'=gY$, the lines in $\PN$ correspond bijectively with the points of the Grassmannian 
$$
\GG(2,N)\,:=\,
\mbox{\raisebox{-.3ex}{$Gl(2)$}$\backslash$\hspace{-.5em} \raisebox{.2ex}{$M(2,N)$}}\,.
$$
Natural homogeneous coordinates on $\GG(2,N)$ are the \emph{Pl\"ucker coordinates}
$Y_{km}$ of the matrix 
$Y=(y_{ij})_{i=1,2;j=1,\ldots,N}$:
\begin{equation}\label{eq:Pluck}
Y_{km}\,:=\,y_{1k}y_{2m}-y_{2k}y_{1m}\,,\qquad k,m=1,\ldots,N\,.
\end{equation}
In more algebraic terms (see \cite{gkz4} p.96), \emph{the ring of homogeneous coordinates on $\GG(2,N)$
in the Pl\"ucker embedding} is the ring
\begin{equation}\label{eq:grasscoord}
\cR_{2,N}\,:=
\modquot{\CC[\Upsilon_{km}\,|\,k,m=1,\ldots,N]\:}{\cI}\,,
\end{equation}
where $\cI$ is the ideal in the polynomial ring $\CC[\Upsilon_{km}]$
generated by the elements
\begin{equation}\label{eq:pluckrel}
\Upsilon_{km}+\Upsilon_{mk}\quad\textrm{and}\quad
\Upsilon_{ij}\Upsilon_{km}+\Upsilon_{ik}\Upsilon_{mj}+\Upsilon_{im}\Upsilon_{jk}
\;,\qquad
\forall i,\,j,\,k,\,m\,.
\end{equation}
$Y_{km}$ is then the element $\Upsilon_{km}\bmod\cI$ of $\cR_{2,N}$.

\

\begin{definition}\label{def:chow form}
\textup{(see \cite{gkz4} pp.99,100,123)}
Let $\cX$ be an irreducible subvariety of $\PN$ of codimension $2$ and degree $d$.
Then the set of lines $\cL_Y$ which intersect $\cX$ is an irreducible hypersurface $\cZ(\cX)$ in the Grassmannian $\GG(2,N)$, called the \emph{associated hypersurface of $\cX$}. 
The hypersurface $\cZ(\cX)$ is given by the vanishing
of a single, up to a multiplicative constant unique, irreducible 
homogeneous element of degree $d$
in $\cR_{2,N}$, the so-called
\emph{Chow form of $\cX$}.
The Chow form of a codimension $2$ subvariety of $\PN$ is the product of the Chow forms of its irreducible components.
\end{definition}

\

Our aim is to explicitly compute the Chow form of the orbit closure
$\ol{\GL[\vu]}$ when $[\vu]\in\PN$ is a point with all homogeneous
coordinates $\neq 0$.

\begin{proposition}\label{proposition:incidence}
For $\vu=(u_1,\ldots,u_N)\in\CNs$ we denote its inverse in $\CNs$ by
$\vu^{-1}=(u^{-1}_1,\ldots,u^{-1}_N)$ 
and its image in $\PN$ by $[\vu]=[u_1:\ldots:u_N]$.
For $Y\in M(2,N)$ we denote by $\vY_{\!\!st}$ the tuple of complex numbers
$\left(Y_{s(e)t(e)}\right)_{e\in\Pe}$, where $Y_{km}$ is the Pl\"ucker coordinate as in (\ref{eq:Pluck}).
Then
\begin{equation}\label{eq:vanish}
\textit{point $[\vu]$ lies on line $\cL_Y$}\quad\Rightarrow\quad
\det\KP(\vY_{\!\!st},\vu^{-1})\,=\,0\,.
\end{equation}
\end{proposition}

\proof
From the fact that a $3\times N$-matrix has rank $\leq 2$ 
if and only if all its $3\times 3$-subdeterminants vanish, we see that
\textit{point $[\vu]$ lies on line $\cL_Y$}
if and only if
\begin{equation}\label{eq:relations1}
u^{-1}_k u^{-1}_m Y_{km}\,+\,
u^{-1}_m u^{-1}_j Y_{mj}\,+\,
u^{-1}_j u^{-1}_k Y_{jk}\,=\,0
\end{equation}
for every triple of elements $j,k,m$ in $\{1,\ldots,N\}$.
Adding equations from the system
(\ref{eq:relations1}) yields an equation
\begin{equation}\label{eq:cycles}
\sum_{i=1}^r
Y_{{k_i}k_{i+1}}u^{-1}_{k_i}u^{-1}_{k_{i+1}}\,=\,0
\end{equation}
for every cycle $(k_1,\ldots,k_r)$, $k_{r+1}=k_1$, of distinct elements of $\{1,\ldots,N\}$. In particular, there is such an equation (\ref{eq:cycles}) 
for the cycle of vertices of a black $2$-cell of $\sP$ and
the terms of that equation match bijectively
with the adjacent white $2$-cells.
This means that in every row of the matrix $\KP(\vY_{\!\!st},\vu^{-1})$
the sum of the entries is $0$.
Consequently $\det\KP(\vY_{\!\!st},\vu^{-1})=0$.
\qed

\

An immediate consequence of Proposition \ref{proposition:incidence}
and Theorem \ref{thm:homogeneity} is:

\begin{corollary}\label{corollary:incidence}
Assume $\LL$ and $\sP$ satisfy Conditions \ref{condition P} and 
\ref{condition LP}.
Let $\vu$ and $Y$ be as in Proposition \ref{proposition:incidence}.
Then $\det\KP(\vY_{\!\!st},\vu^{-1})=0$ if the line $\cL_Y$ intersects the closure $\ol{\GL[\vu]}$ of the orbit $\GL[\vu]$ in $\PN$.
\end{corollary}
\proof
Indeed (\ref{eq:det homogeneity}) and (\ref{eq:vanish})
imply $\det\KP(\vY_{\!\!st},\vu^{-1})=0$ if the line $\cL_Y$ intersects the orbit $\GL[\vu]$. The corollary follows because the set of lines which intersect $\ol{\GL[\vu]}$ is the closure in the Grassmannian $\GG(2,N)$ of the set of lines which intersect $\GL[\vu]$.
\qed

\

In order to remedy for the appearing inverses we now introduce the complementary
bi-adjacency matrix:

\begin{definition}\label{def:complementary}
The \emph{complementary bi-adjacency matrix} $\KP^c(\vz,\vu)$ is
$$
\KP^c(\vz,\vu)\,:=\,
u_1\cdot\ldots\cdot u_N\cdot \KP(\vz,(u^{-1}_1,\ldots, u^{-1}_N)).
$$
with, as usual, $\vu=(u_1,\ldots,u_N)$. The matrix $\KP^c(\vz,\vu)$
has entries in the polynomial ring $\ZZ[\vz,\vu]$.
\end{definition}

\

A further consequence of Proposition \ref{proposition:incidence}
and Theorem \ref{thm:homogeneity} is:

\begin{corollary}\label{corollary:divides}
Assume that $\LL$ and $\sP$ satisfy Conditions \ref{condition P} and \ref{condition LP}.
Let $\vy$ denote the ringhomomorphism (cf. (\ref{eq:grasscoord}))
\begin{equation}\label{eq:pluckerweight}
\vy:\ZZ[z_e\:|\:e\in\Pe\,]\longrightarrow\cR_{2,N}\,,\qquad \vy(z_e)=Y_{s(e)t(e)}\,.
\end{equation}
Then for $\vu\in\CNs$ the 
Chow form of $\ol{\GL[\vu]}$ divides  $\det\KP^c(\vy(\vz),\vu)$ in $\cR_{2,N}$.
\end{corollary}

\proof
There is a slight subtlety because
torsion in $\modquot{\ZN}{\LL}$ makes $\ol{\GL[\vu]}$ reducible. 
If $\cX$ is an irreducible component of $\ol{\GL[\vu]}$ it follows from
Corollary \ref{corollary:incidence} that the Chow form of $\cX$ divides
$\det\KP^c(\vy(\vz),\vu)$, when both are viewed as elements of $\cR_{2,N}$.
Since an irreducible codimension $2$ subvariety of $\PN$ is uniquely determined
by its Chow form (see \cite{gkz4} p.102 Prop.2.5), the Chow forms of different
irreducible components of $\ol{\GL[\vu]}$ are not equal. It follows
that $\det\KP^c(\vy(\vz),\vu)$ is divisible by the Chow form of $\ol{\GL[\vu]}$.
\qed

\begin{theorem}\label{theorem:equal}
Assume that $\LL$ and $\sP$ satisfy Conditions \ref{condition P} and \ref{condition LP}.
Let $\vu\in\CNs$.
\\
Then $\det\KP^c(\vy(\vz),\vu)$ is a 
Chow form for $\ol{\GL[\vu]}$.
\end{theorem}
\textbf{N.B.} We write here \emph{a} Chow form to emphasize that it
is only determined up to a non-zero multiplicative constant.

\proof
In view of Corollary \ref{corollary:divides} we need only prove that the 
Chow form of $\ol{\GL[\vu]}$ and $\det\KP^c(\vy(\vz),\vu)$ are elements
of $\cR_{2,N}$ with the same degree.
\\
Take $\LL^0$ such that
$\LL\subset\LL^0\subset\ZN$ and $\modquot{\LL^0}{\LL}\,=\,(\modquot{\ZN}{\LL})_\tors$.
Set
$$
\GL^0:=\Hom(\modquot{\ZN}{\LL^0},\CC^*)\,,\quad
\va_i^0:= \; \ve_i \bmod\LL^0\,,\quad 
\cA^0:=\{\va_1^0,\ldots,\va_N^0\}\,.
$$
Then $\GL^0$ is the connected component containing the identity in $\GL$.
Via the action of $\CNs$ every irreducible component of
$\ol{\GL[\vu]}$ together with its embedding into $\PN$ is isomorphic to
$\ol{\GL^0\mathsf{1}}\hookrightarrow\PN$, where $\mathsf{1}\in\PN$ has all its components equal to $1$. In \cite{gkz4} p.166 the toric variety $\ol{\GL^0\mathsf{1}}$ is denoted as $X_{\cA^0}$. Since $\LL$ is contained in the kernel
of the map $h$ (see (\ref{eq:h})) the set $\cA^0$ is contained in an
affine hyperplane in $\modquot{\ZN}{\LL^0}\,=\,\ZZ^{N-2}$.
By \cite{gkz4} p.203 Theorem 2.3 and p.99 Proposition 2.2 the degree of the toric variety $X_{\cA^0}\subset\PN$ and the degree of its Chow form
in $\cR_{2,N}$ are both equal to the volume $\mathrm{vol}_{\cA^0}$
of the convex hull of $\cA^0$,
with the volume normalized so that the standard unit cube 
in $\ZZ^{N-3}$ has volume $(N-3)!$ (see \cite{gkz4} p.182).
Let $\va^0_0\in\modquot{\ZN}{\LL^0}$ be the image of
$\va_0\in\modquot{\ZN}{\LL}$. 
It is well-known and easy to prove (see e.g. \cite{S} Eq.(19))
that
$\mathrm{vol}_{\cA^0}\,=\,\frac{1}{2}h(\va^0_0)$.
On the other hand it is obvious that $h(\va_0)h(\va^0_0)^{-1}\,=\,\sharp\left(\modquot{\LL^0}{\LL}\right)$
and that this equals the number of irreducible components of $\ol{\GL[\vu]}$.
Thus the degree of the Chow form of $\ol{\GL[\vu]}$
is equal to $\frac{1}{2}h(\va_0)$. 
It follows from (\ref{eq:bi-adjacency matrix}) and (\ref{eq:sizebw}) that this is also the degree
of $\det\KP^c(\vy(\vz),\vu)$ as an element of $\cR_{2,N}$.
\qed

\section{The principal $\cA$-determinant.}
\label{section:Adeterminant}
In this section we prove Conjecture 10.5 of \cite{S} by combining Theorem \ref{theorem:equal} and  \cite{DS} Proposition 3.2.
The result is:

\begin{theorem}\label{thm:Adeterminant}
Assume that Conditions \ref{condition P} and \ref{condition LP} are satisfied
and that $\modquot{\ZN}{\LL}$ has no torsion. 
Let $B=(b_{ij})_{i=1,2;\,j=1,\ldots,N}$ be a $2\times N$-matrix 
with entries in $\ZZ$ which gives the inclusion
$\LL\hookrightarrow\ZN$ and let 
$\vB_{\!st}$ denote the ring homorphism
$$
\vB_{\!st}:\,\ZZ[z_e\:|\:e\in\Pe]\longrightarrow \ZZ\,,\qquad
\vB_{\!st}(z_e)\,=\,b_{1s(e)}b_{2t(e)}-b_{1t(e)}b_{2s(e)}\,.
$$
Set $\cA=\{\va_1,\ldots,\va_N\}$ and 
let $E_\cA(f)$ denote the
\emph{principal $\cA$-determinant} defined in \cite{gkz4} p.297 Eq. (1.1)
for the Laurent polynomial $f=\sum_{i=1}^Nu_i\vx^{\va_i}$.

Then the following equality holds in the polynomial ring $\ZZ[u_1,\ldots, u_N]$ 
\begin{equation}\label{eq:Adeterminant}
E_\cA(f)\,=\,\pm\vB_{\!st}\left(\det\KP^c(\vz,\vu)\right)\,.
\end{equation} 

\end{theorem}

\proof
One can trace back the definition of the
principal $\cA$-determinant $E_\cA(f)$ from \cite{gkz4} p.297 Eq.(1.1)
to the definition of the Chow form of the subvariety
$X_\cA\,=\,\ol{\GL\mathsf{1}}$ of $\PN$ on \cite{gkz4} p.100.
That is also what Dickenstein and Sturmfels do before defining the
principal $\cA$-determinant, which they call the \emph{full discriminant},
in \cite{DS} Definition 3.1. The result in \cite{DS} Eq. (3.2) is a more or less immediate consequence of this definition. The above Equation (\ref{eq:Adeterminant}) now follows by substituting the expression for the Chow form given in Theorem \ref{theorem:equal} into \cite{DS} Eq. (3.2).
\qed

\section{An example.}\label{section:example}
Figures \ref{fig:fan F1}, \ref{fig:dP3I}
and \ref{fig:bi-ad} show $\LL$, $\sP$, the bi-adjacency matrix
and the verification of Condition \ref{condition LP} for the example
which in the physics literature (e.g. \cite{FHKVW}) is known as \emph{model I of Del Pezzo $3$}. The right-hand picture in Figure \ref{fig:dP3I} is the dual
of the \emph{brane tiling} in Figure 1 of \cite{FHKVW}. One should keep in mind that in \cite{FHKVW} toric data for the singularity are extracted from the \emph{Kasteleyn matrix} (see \cite{FHKVW} Eqs. (3.4), (3.5)),
whereas we work with the bi-adjacency matrix. In this particular example
the Kasteleyn matrix and the bi-adjacency matrix give the same toric data,
but in general they lead to different toric data.

In the context of hypergeometric systems this example is Appell's $F_1$;
see \cite{S} \S 4. The group $\modquot{\ZN}{\LL}$ is torsion free
and can be identified with $\ZZ^4$. This puts $\va_1,\ldots,\va_6$
as points in a hyperplane in $\ZZ^4$.
Figure \ref{fig:classical examples} shows the points $\va_1,\ldots,\va_6$
and their convex hull (often called \emph{the primary polytope})
situated in this $3$-dimensional hyperplane. $\LL$ is the lattice of affine relations between these six points.
\\
From the bi-adjacency matrix in Figure \ref{fig:bi-ad} one easily gets the determinants
\begin{eqnarray*}
\det\KP(\vz,\vu)&=&
\\
&&\hspace{-6em}(z_2z_8z_{12}+z_3z_9z_{10}+z_5z_7z_{11}
-
z_6z_8z_{10}-z_1z_9z_{11}-z_4z_7z_{12})\vu^{[1,1,1,1,1,1]}
\\
&&\hspace{-3em}+\:z_1z_8z_{12}\vu^{[1,2,2,1,0,0]}\:+\:
z_4z_9z_{10}\vu^{[0,0,1,2,2,1]}\:+\:
z_6z_7z_{11}\vu^{[2,1,0,0,1,2]}
\\
&&\hspace{-3em}-\:
z_5z_8z_{10}\vu^{[0,1,2,2,1,0]}\:-\:
z_2z_9z_{11}\vu^{[1,0,0,1,2,2]}\:-\:
z_3z_7z_{12}\vu^{[2,2,1,0,0,1]}\,.
\\[2ex]
\det\KP^c(\vy(\vz),\vu)&=&
\\
&&\hspace{-9em}
(\vY_{56|24|13}+\vY_{12|46|35}+\vY_{34|62|51}
-
\vY_{61|24|35}-\vY_{23|46|51}-\vY_{45|62|13})
\vu^{[2,2,2,2,2,2]}
\\
&&\hspace{-6em}+\:
\vY_{23|24|13}\vu^{[2,1,1,2,3,3]}\:+\:
\vY_{45|46|35}\vu^{[3,3,2,1,1,2]}\:+\:
\vY_{61|62|51}\vu^{[1,2,3,3,2,1]}
\\
&&\hspace{-6em}-\:
\vY_{34|24|35}\vu^{[3,2,1,1,2,3]}\:-\:
\vY_{56|46|51}\vu^{[2,3,3,2,1,1]}\:-\:
\vY_{12|62|13}\vu^{[1,1,2,3,3,2]}\,.
\end{eqnarray*}
with notations $\vu^{[a,b,c,d,e,f]}:=u_1^au_2^bu_3^cu_4^du_5^eu_6^f$
and $\vY_{ab|cd|ef}:=Y_{ab}Y_{cd}Y_{ef}$; here the $Y_{km}$ 
are the Pl\"ucker coordinates on $\GG(2,6)$; see 
(\ref{eq:Pluck})-(\ref{eq:pluckrel}).
These formulas yield for every $\vu\in(\CC^*)^6$ elements
in $\CC[\vz]=\CC[z_1,\ldots,z_{12}]$ and $\cR_{2,6}$, respectively.
Taking the closure of $\{\det\KP(\vz,\vu)\:|\:\vu\in(\CC^*)^6\}$
in $\CC[\vz]$ adds six lines which correspond to the six sides of the secondary polygon in Figure \ref{fig:fan F1}; the top side, for instance, gives
the line
$\{(z_1z_8z_{12}\,x_1-z_5z_8z_{10}\,x_2)\:|\:[x_1:x_2]\in\PP^1\}$.
A similar remark holds for the closure of 
$\{\det\KP^c(\vy(\vz),\vu)\:|\:\vu\in(\CC^*)^6\}$.

\begin{figure}[t]
\begin{center}
\begin{picture}(300,80)(30,10)
\put(140,10){
\begin{picture}(100,80)(0,0)
\put(0,0){\line(1,0){40}}
\put(0,0){\line(0,1){40}}
\put(40,40){\line(-1,0){40}}
\put(40,40){\line(0,-1){40}}
\put(60,20){\line(0,1){40}}
\put(40,0){\line(1,1){20}}
\put(40,40){\line(1,1){20}}
\put(0,40){\line(3,1){60}}
\multiput(0,0)(15,5){4}{\line(3,1){10}}
\put(0,0){\circle*{5}}
\put(40,0){\circle*{5}}
\put(0,40){\circle*{5}}
\put(40,40){\circle*{5}}
\put(60,20){\circle*{5}}
\put(60,60){\circle*{5}}
\put(-15,0){$\va_1$}
\put(45,0){$\va_5$}
\put(-15,40){$\va_4$}
\put(28,30){$\va_2$}
\put(63,20){$\va_3$}
\put(63,60){$\va_6$}
\end{picture}}

\end{picture}
\end{center}
\caption{\label{fig:classical examples}
\textit{\mbox{The set $\cA=\{\va_1,\ldots,\va_6\}$ for 
$\LL=\ZZ^2$ \footnotesize{$\left[\protect\begin{array}{rrrrrr}
\!1&\!1&\!0&\!-1&\!-1&\!0\protect\\ \!0&\!1&\!1&\!0&\!-1&\!-1\protect\end{array}\right]$}}.}}
\end{figure}

On the other hand evaluating $\det\KP^c(\vy(\vz),\vu)$ at the point
$\vl$ amounts in this case to setting all $Y_{km}=1$. This results in
a polynomial in $\CC[u_1,\ldots,u_6]$:
\begin{eqnarray*}
\det\KP^c(\vy(\vz),\vu)_{|\mathrm{at}\:\vl}&=&\\
&&\hspace{-10em}=\,
\vu^{[2,1,1,2,3,3]}+
\vu^{[3,3,2,1,1,2]}+
\vu^{[1,2,3,3,2,1]}-
\vu^{[3,2,1,1,2,3]}-
\vu^{[2,3,3,2,1,1]}-
\vu^{[1,1,2,3,3,2]}
\\
&&\hspace{-10em}=\,
u_1u_2u_3u_4u_5u_6(u_1u_2-u_4u_5)(u_3u_4-u_1u_6)(u_5u_6-u_2u_3).
\end{eqnarray*}
Note that in the first expression for 
$\det\KP^c(\vy(\vz),\vu)_{|\mathrm{at}\:\vl}$
all monomials are invariant for the action of the group $\tGL$ on
$\CC[u_1,\ldots,u_6]$. This invariance does not hold for the factors in the second expression. On the other hand, the factors in this factorization
correspond with the faces of the primary polytope in Figure \ref{fig:classical examples} in agreement with the prime factorization of the principal
$\cA$-determinant in \cite{gkz4} p.299 Theorem 1.2.

There can not be a similar factorization of $\det\KP^c(\vy(\vz),\vu)$
in the ring $\cR_{2,6}[u_1,\ldots,u_6]$, because that would upon specializing
$u_1,\ldots,u_6\in\CC^*$ give a factorization in $\cR_{2,6}$ of
the Chow form of the toric variety 
$\ol{\GL[\vu]}\,\subset\,\PP^5$
and thus contradict the fact that Chow forms of irreducible subvarieties of
$\PP^5$ are irreducible.

\

\noindent\textbf{Acknowledgement.}
I want to express special thanks to Alicia Dickenstein for an e-mail saying that
some results in \cite{S} reminded her of her paper with Bernd Sturmfels
\cite{DS}. This made me aware of the interesting results on Chow forms in \cite{DS,gkz4} and triggered the present paper. Chow forms are functions
on spaces of \emph{algebraic cycles} in projective space. It is therefore
a pleasure to present this work in a volume
in honor of Spencer Bloch.

\end{document}